\documentclass[10pt, titlepage]{amsart}

\usepackage{ae} 
\usepackage[T1]{fontenc}
\usepackage[cp1250]{inputenc}
\usepackage{amsmath}
\usepackage{amssymb, amsfonts,amscd,verbatim}
\usepackage[dvips]{graphicx}
\usepackage{latexsym}
\usepackage{indentfirst}
\usepackage{latexsym}

\usepackage{graphicx}
\usepackage{verbatim}   
\usepackage{color}      
\usepackage{subfigure}  

\usepackage{amsmath}    

\theoremstyle{plain}
\newtheorem{Prop}{Proposition}[section]
\newtheorem{Thm}[Prop]{Theorem}

\newtheorem{Main}[Prop]{Main Theorem}
\newtheorem{ShelahThm}[Prop]{Shelah Theorem}

\theoremstyle{definition}
\newtheorem{Def}[Prop]{Definition}

\theoremstyle{remark}

\newtheorem{Question}[Prop]{\bf Question}

\def\int{\mathop{\roman{int}}}

\def\1{^{-1}}

\def\dokaz{{\bf Proof. }}
\def\edokaz{\hfill $\blacksquare$}

\def\Rips{\mathrm{Rips}}

\numberwithin{equation}{section}

\input pstricks.tex
\input xy
\xyoption{all}

\begin{document}
\title[
Fundamental groups of Peano continua
]%
   {Fundamental groups of Peano continua}
\author{J.~Dydak}
\address{University of Tennessee, Knoxville, TN 37996, USA}
\email{dydak@math.utk.edu}

\author{\v Z.~Virk}
\address{University of Tennessee, Knoxville, TN 37996, USA}
\email{zigavirk@gmail.com}

\date{ \today
}
\keywords{coarse geometry, coarse connectivity, finitely presented groups, fundamental group, locally connected compact metric spaces}

\subjclass[2000]{Primary 55Q52; Secondary 20F65, 14F35, 14B0}

\thanks{Supported in part by the Slovenian-USA research grant BI--US/05-06/002 and the ARRS
research project No. J1--6128--0101--04}
\thanks{The first-named author was partially supported
by MEC, MTM2006-0825.}

\begin{abstract}

Extending a theorem of Shelah we prove that fundamental groups of Peano continua
(locally connected and connected metric compact spaces) are finitely presented if they are countable.
The proof uses ideas from geometric group theory.

\end{abstract}

\maketitle

\medskip
\medskip
\tableofcontents
\section{Introduction}

This paper is motivated by a question posed to the second author by Mladen Bestvina
during his talk at the Spring Topology and Dynamics Conference in Gainesville (March 7-9, 2009):

\begin{Question}\label{MainQuestion}
Is the fundamental group of a Peano continuum finitely presented
if it is countable?
\end{Question}

It turns out that question was also posed by 
de la Harpe \cite{Har} on p.48 and it is relevant in view of the following:

\begin{ShelahThm}\label{ShThm}
If $X$ is a Peano continuum and $\pi_1(X)$ is countable, then $\pi_1(X)$ is finitely generated.
\end{ShelahThm}

Pawlikowski \cite{Paw} presented another proof of \ref{ShThm}
from which we extract the following (see the paragraph preceding Lemma 2
in \cite{Paw} or Theorems 2 and 8 in \cite{Fab}):

\begin{Thm}[Pawlikowski \cite{Paw}]\label{PawThm}
If $X$ is a Peano continuum and $\pi_1(X)$ is countable, then $X$ is semi-locally simply connected.
\end{Thm}

Notice that the second author constructed (see \cite{Vir}), for each countable group $G$, 
a $2$-dimensional path-connected subcontinuum $X_G$ of $R^4$ whose fundamental group
is $G$ (see \cite{KR} and \cite{Prze} for earlier constructions of compact spaces
with a given fundamental group).

Our solution to \ref{MainQuestion} is based on an application of methods of 
geometric group theory: we construct a geometric action of $\pi_1(X)$ on a coarsely $1$-connected proper geodesic space $\tilde X$ and we use \v Svarc-Milnor Lemma
(\cite[page 140]{BH}) plus
the fact $G$ is finitely presented if and only if it is coarsely $1$-connected.

In geometric group theory, a geometry is any proper, geodesic metric space. An action of a finitely-generated group $G$ on a geometry $X$ is {\bf geometric} if it satisfies the following conditions:
\begin{enumerate}
\item Each element of $G$ acts as an isometry of $X$.
\item The action is cocompact, i.e. the quotient space $X/G$ is a compact space.
\item The action is properly discontinuous, with each point having a finite stabilizer.
\end{enumerate}

Let a group $G$ act on a topological space $X$ by homeomorphisms. 
Consider a subgroup $H \subset G$. One then says that a set $Y$ is precisely invariant under $H$ in $G$ if

 $$   \forall h \in H, \quad h(Y)=Y \;\mbox{ and }\; \forall g \in G-H, \quad gY \cap Y = \varnothing.$$

Then let $G_x$ be the stabilizer of $x$ in $G$. One says that $G$ acts {\bf discontinuously} at $x$ in $X$ if the stabilizer $G_x$ is finite and there exists a neighborhood $U$ of $x$ that is precisely invariant under $G_x$ in $G$. If $G$ acts discontinuously at every point $x$ in $X$, then one says that $G$ acts {\bf properly discontinuously} on $X$.

\begin{Thm}[\v Svarc-Milnor \cite{BH} or \cite{BDM1}] \label{Svarc-Milnor}
A group $G$ acting properly discontinuously and cocompactly via isometries on a
length space $X$ is finitely generated and induces a
quasi-isometry equivalence $g\to g\cdot x_0$ for any $x_0\in X$.
\end{Thm}

{\bf Added in proof:}
We were informed by Greg Conner \cite{CPC} that Katsuya Eda answered \ref{MainQuestion} about 5 years ago
(unpublished).  The
argument is contained in \cite{CC}: by \ref{PawThm}  the space
is semi-locally simply-connected and homotopically Hausdorff. Corollary 5.7 says such
space has finitely presented fundamental group.

Alternatively, it is pointed out in Lemma A.3 in \cite{CL} that this
shows that such a space has finitely presented fundamental group.
Here is the argument from Lemma A.3.  Semilocally simply-connected
implies that the space is two-set simple -- see \cite{CA}.  This implies
that the fundamental group is group is the fundamental group of the
nerve of a finite cover which implies that it's finitely presented,
again from \cite{CA}.

\section{Coarse 1-connectivity of uniformly path connected spaces}

In order to complete the proof of our main result \ref{MainThm}, we need to relate coarse $1-$connectedness to simple connectedness.

Recall $(X,d)$ is {\bf coarsely $1$-connected} if it is $t$-chain connected
for some $t > 0$ (that is equivalent to $(X,d)$ being coarsely $0$-connected - see \cite{Kap}) and for each $r > 0$ there is $R > 0$ such that the induced map
$\Rips_r(X)\to \Rips_R(X)$ induces the trivial homomorphism of the fundamental groups
(see Definition 42 on p.19 of \cite{Kap}).
Here $\Rips_r(X)$ is the {\bf Rips complex} of $X$, i.e. the complex whose simplices are all finite subsets $A$ of $X$ of diameter at most $r$.

\begin{Def}
A path connected metric space $(X,d)$ is \textbf{uniformly path connected } if there is a function $\alpha\colon (0,\infty)\to (0,\infty)$ so that every two points $x,y\in X$  can be connected by a path of diameter at most $\alpha (d(x,y))$.

The fundamental group $\pi_1(X,x_0)$ of a path connected metric space $X$ is \textbf{uniformly generated} (see \cite{FW}) if it  has a generating set of loops of diameter at most $R$ for some $R>0.$ Equivalently,  every map $f\colon (S^1,0)\to (X,x_0)$ can be extended  over $1-$skeleton of some subdivision $\tau$ of $(B^2,0)$ to a map $F$ so that the diameter $F(\partial \Delta)$ is at most $R$, for every simplex $\Delta$ of $\tau$.
\end{Def}

\begin{Thm}\label{Coarse1Connectivity}
Suppose $X$ is a uniformly path connected space. $X$ is coarsely $1-$connected
if and only if $\pi_1(X,x_0)$ is uniformly generated.
\end{Thm}
\dokaz Assume $X$ is coarsely $1-$connected.
Fix positive numbers $r,R$ so that $\pi_1$ applied to $\Rips_r(X)\to \Rips_R(X)$ is trivial. Furthermore, let $l$ be a positive number so that every two points of $X$ that are at most $R$ apart can be connected by a path of diameter at most $l$. Let $\alpha \colon (S^1,0)\to (X,x_0)$ be a loop. Subdivide $(S^1,0)$ to obtain a subdivision $\tau$ (notation: $S^1_\tau$) so that the diameter of $\alpha(\Delta)$ is at most $r$ for every edge $\Delta$ of $\tau$. The map $\alpha|_{(S^1_\tau)^{(0)}}$ induces a simplicial map $\tilde \alpha \colon (S^1_\tau,0)\to (\Rips_r(X),x_0)$, which extends to a map $\tilde \beta\colon (B^2,0)\to \Rips_R(X)$. We may assume $\sigma$ is a subdivision of $(B^2,0)$ so that $\tilde \beta$ is simplicial and $\sigma|_{S^1}$ is a subdivision of $\tau.$  Then $\tilde \beta$ induces a map $\beta \colon ((B^2_{\sigma})^{(1)},0)\to (X,x_0)$ as follows: 
$\beta$ equals $\tilde\beta$ on vertices and $S^1$ and for every edge $E$ of $B^2_\sigma\setminus S^1$ we can connect  two boundary points $\beta(\partial E)$ by a path of diameter at most $l$. Hence we obtain an extension $\beta \colon ((B^2_{\sigma})^{(1)},0)\to (X,x_0)$ of $\alpha$ so that diameter of $\beta(\Delta)$ is at most $2\cdot l$  for every simplex $\Delta$ of $B^2_{\sigma}$. This means that $\pi_1(X,x_0)$ is $2\cdot l$-generated.

Assume $\pi_1(X,x_0)$ is uniformly generated by loops of diameter at most $D$. Fix $r> 0, l>0$ so that every two points of distance at most $r$ can be connected by a path of diameter at most $l$. We can assume $D>l$. Pick any simplicial map $\alpha \colon (S^1_\tau,0)\to (\Rips_r(X),x_0)$. It induces a map $\tilde \alpha \colon ((S^1_\tau)^{(0)},0)\to (X,x_0)$ as follows: For every edge $E$ of $S^1_\tau$ we connect  two boundary points $\tilde \alpha(\partial E)$ by a path of diameter at most $l$ to obtain a map $\tilde \alpha \colon (S^1_\tau,0)\to (X,x_0)$. Such map extends over $1-$skeleton of some subdivision $\sigma$ (containing $\tau$) of $(B^2,0)$ to a map $\tilde \beta$ so that diameter $\tilde \beta(\partial \Delta)$ is at most $D$, for every simplex $\Delta$ of $\sigma$. Then $\tilde \beta $ induces a map $\beta \colon ((B^2_\sigma)^{(0)},0)\to \Rips_{D}(X)$ which extends over $B^2_\sigma$. Note that $\beta|_{(\partial B^2,0)} \simeq \alpha$: for every edge $E$ of $\tau$ the set $\beta(E)\cup \alpha(E)$ is contained in a simplex of $\Rips_{D}(X)$ because of uniform path connectedness and $D>l.$
\edokaz

\section{Main result}

Given a Peano continuum $X$ we assume it has a geodesic metric $d_X$ (see \cite{Bing}). Pick a base point $x_0$ of $X$ and consider the space $\tilde X$
of homotopy (rel.endpoints) classes of paths in $X$ originating at $x_0$.

In this section we assume $X$ is semi-locally simply connected.

\begin{Def}
Given $[\alpha]\in\tilde X$ and a path $\beta$ in $X$ originating at $\alpha(1)$,
the {\bf canonical lift} $\tilde\beta$ of $\beta$ is a path in $\tilde X$
defined by $\tilde \beta(t)=[\alpha\ast (\beta\vert_{[0,t]})]$, the concatenation of $\alpha$
and $\beta$ restricted to interval $[0,t]$.
\end{Def}

Given two elements $[\alpha]$ and $[\beta]$ of $\tilde X$ we define the distance
$d([\alpha],[\beta])$ as the infimum of lengths $l(\gamma)$ of all paths
$\gamma$ from $\alpha(1)$ to $\beta(1)$ such that $\gamma$
is homotopic rel.endpoints to $\alpha^{-1}\ast\beta$.

\begin{Prop}\label{dIsAMetric} $(\tilde X,d)$ is a proper geodesic space
such that the endpoint projection $p\colon \tilde X\to X$
is $1$-Lipschitz and canonical lifts of geodesics in $X$ are geodesics in $\tilde X$.
\end{Prop}
\dokaz Let $\delta > 0$ be a number such that any loop in $X$ of diameter less than $4\cdot \delta$ is null-homotopic in $X$. Notice that any two paths at distance less than $\delta$ are homotopic rel.endpoints if they join the same two points.

Given two elements $[\alpha], [\beta]$ of $\tilde X$ the path $\alpha^{-1}\ast \beta$
can be approximated by a piecewise-geodesic path $\gamma$. As $l(\gamma)$ is finite,
so is $d([\alpha],[\beta])$. If $d([\alpha],[\beta])=0$, then $\alpha(1)=\beta(1)$. As $d([\alpha],[\beta])=0$ there is a loop
$\gamma$ at $x_1$ of length less than $\delta$ satisfying $\gamma\sim \alpha^{-1}\ast\beta$.
That means $\alpha\sim\beta$ as $\gamma$ is null-homotopic in $X$.
Thus $[\alpha]=[\beta]$ if $d([\alpha],[\beta])=0$.
It is easy to see $d$ is symmetric and satisfies the Triangle Inequality.
\par Notice $d([\alpha],[\beta])\ge d_X(\alpha(1),\beta(1))$, so $p$ is $1$-Lipschitz.
Also, it is clear that canonical lifts of geodesics in $X$ are geodesics in $\tilde X$.
\par Suppose $\gamma_n$ is a sequence of paths in $X$ joining $\alpha(1)$
and $\beta(1)$ such that $l(\gamma_n)$ converges to $M=d([\alpha],[\beta])$
and $\gamma_n\sim\alpha^{-1}\ast\beta$ for all $n\ge 1$.
We may assume each $\gamma_n$ is parametrized so that the length of $\gamma_n\vert_{[0,t]}$ is $t\cdot l(\gamma_n)$.
Subdivide the interval $[0,1]$ into points $y_0=0,y_1,\ldots,y_k=1$ such that $0 < y_{i+1}-y_i <  \frac{\delta}{2\cdot M}$ for all $0\leq i < k$.
We may assume $\gamma_n(y_i)$ converges to $z_i\in X$ for each $0\leq i\leq k$.
The piecewise-geodesic path $\omega$ from $\alpha(1)$ to $\beta(1)$
obtained by connecting points $z_0$, $z_1$, \ldots,$z_k$ is homotopic to $\gamma_n$
for $n$ large enough. Also, $l(\omega)$ equals the limit of $l(\gamma_n)$,
so $l(\omega)=d([\alpha],[\beta])$. Notice the canonical lift of $\omega$ is a geodesic
from $[\alpha]$ to $[\beta]$ in $\tilde X$.
\par To show $(\tilde X,d)$ is a proper metric space assume
$\{[\alpha_n]\}_{n\ge 1}$ is a bounded sequence in $\tilde X$.
We may assume $\alpha_n(1)$ converges to $x_1$ and then alter each $\alpha_n$ by concatenating it with the geodesic from $\alpha_n(1)$ to $x_1$.
It suffices to show that the resulting sequence of elements $[\beta_n]$ of $\tilde X$
has a convergent subsequence.
First of all, we may assume the sequence of lengths $l(\beta_n)$ converges to $M > 0$
(if $M=0$, then $[\beta_n]$ converge to $[c]$), each $\beta_n$ is piecewise-geodesic
and the length of $\beta_n\vert_{[0,t]}$ is $t\cdot l(\beta_n)$. 
Subdivide the interval $[0,1]$ into points $y_0=0,y_1,\ldots,y_k=1$ such that $0 < y_{i+1}-y_i < \frac{\delta}{2\cdot M}$ for all $0\leq i < k$.
We may assume $\beta_n(y_i)$ converges to $z_i\in X$ for each $0\leq i\leq k$.
The piecewise-geodesic path $\omega$ from $x_0$ to $x_1$
obtained by connecting points $z_0$, $z_1$, \ldots,$z_k$ is homotopic to $\beta_n$
for $n$ large enough. That means $[\beta_n]$ is constant starting from a sufficiently large $n$.
\edokaz

\begin{Prop}\label{TildeXIsGit} $(\tilde X,d)$ is a simply connected and the endpoint projection $p\colon \tilde X\to X$
is a covering map.
\end{Prop}
\dokaz Let $\delta > 0$ be a number such that any loop in $X$ of diameter less than $4\cdot \delta$ is null-homotopic in $X$.

\par {\bf Claim:} For any $[\alpha]\in\tilde X$ the restriction of $p$
to the ball $B([\alpha],\delta)$ is an isometry onto $B(\alpha(1),\delta)$.
\par {\bf Proof of Claim:} Given $\beta,\omega\in B([\alpha],\delta)$
let $\gamma$ be the geodesic path from $\beta(1)$ to $\omega(1)$.
As $d([\beta],[\omega]) < 2\cdot\delta$ there is a path $\lambda$ from
$\beta(1)$ to $\omega(1)$ of length less than $2\cdot\delta$
such that $\lambda\sim\beta^{-1}\ast\omega$. 
Observe $\lambda\sim\gamma$ as both paths are of diameter less than $2\cdot\delta$.
That means $d([\beta],[\omega])=d_X(\beta(1),\omega(1))$
as the length of $\gamma$ equals $d_X(\beta(1),\omega(1))$
and $d([\beta],[\omega])\ge d_X(\beta(1),\omega(1))$.
\par Given $[\beta]\in \tilde X$ with $\beta(1)\in B(x_1,\delta)$
let $\gamma$ be the geodesic path from $\beta(1)$ to $x_1$.
Observe $d([\beta],[\beta\ast\gamma]) < \delta$
and $p([\beta\ast\gamma])=x_1$. That means $p^{-1}(B(x_1,\delta))$
is the union of balls $B([\alpha],\delta)$ with $\alpha$ ranging over
all paths in $p^{-1}(x_1)$. By Claim we conclude $p$ is a covering projection.
\par To show $\tilde X$ is simply connected suppose $\alpha$ is a loop
in $\tilde X$ based at the trivial path. Since $p(\alpha)$ can be homotoped
to a piecewise-geodesic loop and canonical lifts of piecewise-geodesic loops are paths in $\tilde X$,
we may assume $\alpha$ is the canonical lift of a piecewise-geodesic loop
$\beta$ based at $x_0$. The canonical lift of $\beta$ is a loop if and only if $\beta$ is null-homotopic.
As $p$ is a covering projection, $\alpha$ is null-homotopic as well.
\edokaz

\begin{Prop}\label{ActionIsGood}
The action of $G=\pi_1(X,x_0)$ on $\tilde X$ ($g\cdot [\alpha]$ being
$[\beta\ast\alpha]$, where $[\beta]=g$) is geometric.
\end{Prop}
\dokaz $G$ acts by isometries as $d(g\cdot [\alpha],g\cdot [\beta])=
d([\alpha],[\beta])$ for all $\alpha,\beta\in\tilde X$.
\par Let $\delta > 0$ be a number such that any loop in $X$ of diameter less than $4\cdot \delta$ is null-homotopic in $X$.
Given $[\alpha]\in\tilde X$ let $U$ be the $\delta$-ball around $[\alpha]$ in $\tilde X$.
If $[\beta]\in U\cap (g\cdot U)$ there are paths $\gamma_i$, $i=1,2$,
such that $\beta\sim \alpha\ast\gamma_1$,
$\beta\sim g\cdot \alpha\ast\gamma_2$ and $l(\gamma_i) < \delta$ for $i=1,2$.
Thus $g=[\alpha\ast \gamma_1\ast \gamma_2^{-1}\ast\alpha^{-1}]$ equals $1$ in $G$
as $\gamma_1\sim\gamma_2$ (both are paths of diameter less than $\delta$
and join the same points). That proves the action of $G$ on $\tilde X$ is properly discontinuous.
\par Since $p\colon\tilde X\to X$ is open and, set-theoretically, equals
$\tilde X\to \tilde X/G$, $\tilde X/G$ is homeomorphic to $X$ proving that the action of $G$ on $\tilde X$ is cocompact.
\edokaz

\begin{Main}\label{MainThm}
The fundamental group of a Peano continuum $X$ is finitely presented if it is countable.
\end{Main}
\dokaz By the \v Svarc-Milnor Lemma and \ref{ActionIsGood} the group $G=\pi_1(X,x_0)$
is finitely generated and
is quasi-isometric to $\tilde X$. As $\tilde X$ is coarsely 1-connected (see \ref{Coarse1Connectivity})
and coarse 1-connectivity is an invariant of quasi-isometries (see Corollary 47 in \cite{Kap}), $G$ is also coarsely 1-connected.
As $G$ is a finitely generated group it means $G$ is finitely presented (see the proof of Corollary 51 in \cite{Kap} on p.22 or Proposition 8.24 in \cite{BH}).
Alternatively, the fundamental group of the Cayley graph $\Gamma(G)$ of $G$ must be uniformly generated by \ref{Coarse1Connectivity} which means $G$ is finitely presented.
\edokaz

\end{document}